\documentclass[final,leqno,letterpaper]{etna}

\usepackage{graphicx,amsmath,amsfonts,amssymb}

\usepackage{bm}
\usepackage{xcolor}

\makeatletter
\def\BState{\State\hskip-\ALG@thistlm}
\makeatother

\newcommand{\rev}[2][black]{{\textcolor{#1}{#2}}}
\newcommand*{\de}{\mathop{}\!\mathrm{d}}

\def\ba{\begin{aligned}}

\def\be{\begin{equation}}

\def\cF{{\mathcal\FFF}}
\def\cl {\nonumber \\}
\def\cN{{\mathcal\NNN}}

\def\ea{\end{aligned}}

\def\ee{\end{equation}}

\def\el {\nonumber }

\def\FFF{{F}}
\def\g{{\bm g}}

\def\kkk{{k}}
\def\KKK{{K}}

\def\LLL{{L}}
\def\LLLk{{\LLL_\kkk}}
\def\MMM{{M}}
\def\n{{\bm n}}

\def\NNN{{N}}

\def\R{{\mathbb R}}

\def\uuu{{u}}

\def\uuul{{\uuu_\LLL}}
\def\uuulk{{\uuu_\LLLk}}
\def\uuuh{{\uuu_\NNN}}

\def\vvv{{v}}
\def\VVV{{V}}
\def\vvvl{{\vvv_\LLL}}
\def\VVVl{{\VVV_\LLL}}
\def\VVVlk{{\VVV_\LLLk}}

\def\VVVh{{\VVV_\NNN}}

\newcommand{\mmu}{\boldsymbol{\mu}}

\newtheorem{remark}{Remark}

\setbibdata{1}{xx}{46}{2017} 

\hypersetup{%
    pdftitle={Document title},
    pdfauthor={John Doe, Erwin Schr\"odinger},
    pdfkeywords={a keyword, another keyword, and another one}
    }

\title{A Comparison of Reduced-Order Modeling Approaches \rev{Using Artificial Neural Networks} for PDEs with Bifurcating Solutions \thanks{%
Received... Accepted... Published online on... Recommended by.... 
We acknowledge the support by European Union Funding for Research and Innovation - Horizon 2020 Program - in theframework of European Research Council Executive Agency: Consolidator Grant H2020 ERC CoG 2015 AROMA-CFD project 681447 “Advanced Reduced Order Methods with Applications in Computational Fluid Dynamics” (P.I.Prof. Gianluigi Rozza).
We  also  acknowledge  the  INDAM-GNCS  project  “Advanced  intrusive  and  non-intrusive  model  order  reduction techniques and applications”, 2019.
 We acknowledge project MIUR PRIN 2017 NA-FROM-PDEs.
This work was also partially supported by US National Science Foundation through grant DMS-1620384 and DMS-195353.
}}

\author{Martin W. Hess\footnotemark[2]
        \and Annalisa Quaini\footnotemark[3] \and Gianluigi Rozza\footnotemark[4] }

\shorttitle{A Comparison of ROM Approaches for PDEs with Bifurcating Solutions} 
\shortauthor{M. Hess, A. Quaini, G. Rozza}

\begin{document}


\maketitle

\renewcommand{\thefootnote}{\fnsymbol{footnote}}

\footnotetext[2]{SISSA Mathematics Area, mathLab, International School for Advanced Studies, via Bonomea 265, I-34136 Trieste, Italy.}
\footnotetext[3]{Department of Mathematics, University of Houston, Houston, Texas 77204, USA.}
\footnotetext[4]{SISSA Mathematics Area, mathLab, International School for Advanced Studies, via Bonomea 265, I-34136 Trieste, Italy.}

\begin{abstract}
This paper focuses on reduced-order models (ROMs) built for the efficient treatment of 
PDEs having solutions that bifurcate as the values of multiple input parameters change. 
First, we consider a method called local ROM that uses k-means algorithm to cluster snapshots
and construct local POD bases, one for each cluster. We investigate one key ingredient of this approach:
the local basis selection criterion. Several criteria are compared and it is found that a criterion 
based on a regression artificial neural network (ANN) provides the most accurate results
for a channel flow problem exhibiting a supercritical pitchfork bifurcation.
The same benchmark test is then used to compare the local ROM approach with the
regression ANN selection criterion to an established global projection-based ROM and 
a recently proposed ANN based method called POD-NN. We show that our local ROM
approach gains more than an order of magnitude in accuracy over the global projection-based ROM. 
However, the POD-NN provides consistently more accurate approximations than the local projection-based ROM.
\end{abstract}

%
%

%




\section{Introduction}

We consider the problem of finding a function $\uuu\in\VVV$ such that 
\be\label{abs-pro}
    \cN(\uuu,\vvv;\mmu) = \cF(\vvv;\mmu)\quad\forall\, \vvv\in\VVV, 
\ee
where $\mmu\in D$ denotes a point in a parameter domain $D\subset \R^\MMM$, $\VVV$ a function space, 
$\cN(\cdot,\cdot;\mmu)$ a given form that is linear in $\vvv$ but generally nonlinear in $\uuu$, and $\cF(\cdot)$ 
a linear functional on $\VVV$. Note that either $\cN$ or $\cF$ or both could depend on some or all the 
components of the parameter vector $\mmu$. We view the problem \eqref{abs-pro} as a variational 
formulation of a nonlinear partial differential equation (PDE) or a system of such equations 
in which $\MMM$ parameters appear. 
We are interested in systems that undergo bifurcations,i.e.~the solution $\uuu$ of \eqref{abs-pro} 
differs in character 
for parameter vectors $\mmu$ belonging to different subregions of $D$. 
We are particularly interested in situations that require solutions of \eqref{abs-pro} for a set of parameter vectors that span across 
two or more of the subregions of the bifurcation diagram. This is the case, for example, if
one needs to trace the bifurcation diagram. 

In general, one could approximate the solutions to \eqref{abs-pro}
using a Full Order Method (FOM), like for example the Finite Element Method
or the Spectral Element Method. 
Let $\VVVh$ be a $\NNN$-dimensional subspace that is a subset of $\VVV$.
A FOM seeks an approximation $\uuuh\in\VVVh$ such that
\be\label{abs-dis}
    \cN_h(\uuuh,\vvv;\mmu) = \cF_h(\vvv;\mmu)\quad\forall\, \vvv\in\VVVh, 
\ee
where $\cN_h$ and $\cF_h$ are the discretized forms for $\cN$ and $\cF$.
FOMs are often expensive, especially if multiple solutions are needed. 
For this reason, one is interested in finding surrogate methods that are much less costly.
Such surrogates are constructed using a ``few'' solutions obtained with the FOM. 
Here, we are interested in {\em reduced-order models} (ROMs) for which one constructs a
 low-dimensional approximating subspace $\VVVl\subset\VVVh$ of dimension $\LLL$ 
that still contains an acceptably accurate approximation $\uuul$ to the \rev{FOM solution $\uuuh$, and thus also to the} solution $u$ of \eqref{abs-pro}. 
That approximation is determined from the reduced discrete system
\be\label{abs-romg}
    \cN_h(\uuul,\vvv;\mmu) = \cF_h(\vvv;\mmu)\quad\forall\, \vvvl\in\VVVl 
\ee
that, if $\LLL\ll\NNN$, is much cheaper to solve compared to \eqref{abs-dis}. 
We refer to approach \eqref{abs-romg} as \emph{global ROM}, since
a single {\em global} basis is used to determine the ROM approximation $\uuul$ 
at any chosen parameter point $\mmu\in D$.

Global ROMs in the setting of bifurcating solutions are considered in the early papers 
\cite{NOOR1982955,Noor:1994,Noor:1983,NOOR198367} for buckling bifurcations in solid mechanics. 
More recently, in \cite{Terragni:2012} it is shown that a 
Proper Orthogonal Decomposition (POD) approach allows for considerable computational time 
savings for the analysis of bifurcations in some nonlinear dissipative systems. 
Reduced Basis (RB) methods have been used to study symmetry breaking bifurcations
\cite{Maday:RB2,PLA2015162} and Hopf bifurcations  \cite{PR15} for
natural convection problems. A RB method for symmetry breaking bifurcations in contraction-expansion channels
has been proposed in  \cite{PITTON2017534}. 
A RB method for the stability of flows under perturbations in the forcing term or in the boundary conditions,
is introduced  in \cite{Yano:2013}. Furthermore, in~\cite{Yano:2013} 
it is shown how a space-time inf-sup constant approaches zero as the computed solutions get close to a bifurcating value. 
Recent works have proposed ROMs for bifurcating solutions in structural mechanics \cite{pichirozza}
and for a nonlinear Schr\"{o}dinger equation, called Gross--Pitaevskii equation \cite{Pichi2020}, respectively.
Finally, we would like to mention that machine learning techniques based on sparse optimization 
have been applied to detect bifurcating
branches of solutions for a two-dimensional laterally heated cavity and Ginzburg-Landau model
in \cite{BTBK14,KGBNB17}, respectively.
Finding all branches after a bifurcation occurs can be done with deflation methods (see \cite{PintorePichiHessRozzaCanuto2019}), which require
introducing a pole at each known solution. In principle, the ROM approaches under investigation here can be combined with deflation techniques.

In a setting of a bifurcation problem (i.e.~$D$ consists of subregions for which the corresponding solutions of \eqref{abs-pro} 
have different character), it may be the case that $\LLL$, 
although small compared to $\NNN$, may be large enough so that solving system 
\eqref{abs-romg} many times becomes expensive.
To overcome this problem, in \cite{Hess2019CMAME} we proposed a
\emph{local ROM} approach. The idea is to construct several {\em local} bases
(in the sense that they use solutions for parameters that lie in subregions of the parameter domain), each of which is 
used for parameters belonging to a different subregion of the bifurcation diagram. 
So, we construct $\KKK$ such local bases of dimension $\LLLk$, each spanning a local subspace $\VVVlk\subset\VVVh$. 
We then construct $\KKK$ {\em local reduced-order models} 
\be\label{abs-roml}
\ba
    \cN(\uuulk,\vvv;\mmu) = \cF(\vvv;\mmu)\quad\forall\, \vvv\in\VVVlk &\qquad \mbox{for $\kkk=1,\ldots,\KKK$}
\ea
\ee
that provide  acceptably accurate approximations $\uuulk$ to the solution $u$ of \eqref{abs-pro} 
for parameters $\mmu$ belonging to different parts of the bifurcation diagram. 
A key ingredient in this approach is how to identify which local basis should be used in \eqref{abs-roml} 
to determine the corresponding ROM approximation for any parameter point $\mmu\in D$ 
that was not among those used to generate the snapshots.
Several criteria are proposed and compared for a 
two-parameter study in Sec.~\ref{sec:criteria}. 
To the best of our knowledge, no work other than \cite{Hess2019CMAME} addresses
the use of local ROM basis for bifurcation problems. This is the continuation of previous work on model reduction with spectral element methods \cite{HessRozza2019}
and including parametric variations of the geometry \cite{10.1007/978-3-030-39647-3_45}, \cite{HessQuainiRozza2020}.

This paper aims at comparing one global ROM approach, 
our local ROM approach with the ``best'' criterion to select 
the local basis, and a recently proposed RB method that uses 
neural networks to accurately approximate the coefficients of the reduced model
\cite{HESTHAVEN201855}. This third method is referred to as POD-NN.
The global ROM as explained in, e.g., uses the most dominant POD modes of a uniform
snapshot set over the parameter domain. 
The dominant modes define the projection space for every parameter evaluation of interest and is in this sense \emph{global}
with respect to the parameter space.
See, e.g., \cite{LMQR:2014} for more details. Like the global ROM, the
POD-NN employs the most dominant POD modes. However, the difference is that  
 the coefficients of the snapshots in the ROM expansion are
used as training data for an artificial neural network.
\rev{The local ROM first employs a classification ANN to determine the corresponding cluster of a parameter location, but this can be improved upon with a regression ANN using the relative errors of the local ROMs at the snapshot locations as training data.}
As a concrete setting for the comparison, we use the Navier-Stokes equations
and in particular flow through a channel with a contraction. 

The outline of the paper is as follows. In Sec.~\ref{sec:NS}, we briefly
present the Navier-Stokes equations and consider a specific benchmark test. 
Sec.~\ref{sec:num_res} reports the comparison of many local basis selection
criteria for the local ROM approach and the main comparison of the three
ROM approaches. Concluding remarks are provided in Sec.~\ref{sec:concl}.

\section{Application to the incompressible Navier-Stokes equations}\label{sec:NS}

The Navier-Stokes equations describe the incompressible motion of a viscous, Newtonian fluid in the spatial domain $\Omega \subset\mathbb{R}^d$, $d = 2$ or $3$, over a time interval of interest $(0, T]$. 
They are given by
\begin{equation}\label{NS-1}
\begin{aligned}
\frac{\partial {\bm u}}{\partial t}+({\bm u}\cdot  \nabla {\bm u})-\nu\Delta {\bm u}+\nabla p={\bm 0} &\qquad\mbox{in } \Omega\times(0,T]\\
\nabla \cdot {\bm u}=0&\qquad\mbox{in } \Omega\times (0,T],
\end{aligned}
\end{equation}
where ${\bm u}$ and $p$ denote the unknown velocity and pressure fields, respectively, and
$\nu>0$ denotes the kinematic viscosity of the fluid. Note that there is no external body force
because we will not need one for the specific benchmark test under consideration. 

Problem \eqref{NS-1} needs to be endowed with initial and boundary conditions, e.g.:
\begin{eqnarray}
{\bm u}={\bm u}_0&\qquad&\mbox{in } \Omega \times \{0\} \label{IC}\\
{\bm u}={\bm u}_D&\qquad&\mbox{on } \partial\Omega_D\times (0,T] \label{BC-D} \\
-p \n + \nu \frac{\partial {\bm u}}{\partial \n}=\g &\qquad&\mbox{on } \partial\Omega_N\times, (0,T], \label{BC-N}
\end{eqnarray}
where $\partial\Omega_D \cap \partial\Omega_N = \emptyset$ and $\overline{\partial\Omega_D} \cup \overline{\partial\Omega_N} =  \overline{\partial\Omega}$. Here, ${\bm u}_0$, ${\bm u}_D$, and $\g$ are given and $\n$ denotes the unit normal vector on the boundary $\partial\Omega_N$ directed outwards. In the rest of this section, we will explicitly denote the dependence of the solution of the problem \eqref{NS-1}-\eqref{BC-N} on the parameter vector $\mmu$.

Let $L^2(\Omega)$ denote the space of square integrable functions in $\Omega$ and $H^1(\Omega)$ the space of functions belonging to $L^2(\Omega)$ with first derivatives in $L^2(\Omega)$. Moreover, let 
\begin{align}
{\bm V} &:= \left\{ {\bm v} \in [H^1(\Omega)]^d: ~ {\bm v} = {\bm u}_D \mbox{ on }\partial\Omega_D \right\}, \cl
{\bm V}_0 &:=\left\{{\bm v} \in [H^1(\Omega)]^d: ~ {\bm v} = \boldsymbol{0} \mbox{ on }\partial\Omega_D \right\}. \el
\end{align}
The standard variational form corresponding to \eqref{NS-1}-\eqref{BC-N} is: 
find $({\bm u}(\mmu),p(\mmu))\in {\bm V} \times L^2(\Omega)$ satisfying the initial condition \eqref{IC} and
\begin{equation}\label{eq:weakNS-1}
\begin{aligned}
&\int_{\Omega} \frac{\partial{\bm u}(\mmu)}{\partial t}\cdot{\bm v}\de\mathbf{x}+\int_{\Omega}\left({\bm u}(\mmu)\cdot\nabla {\bm u}\right)\cdot{\bm v}\de \mathbf{x} - \int_{\Omega}p(\mmu)\nabla \cdot{\bm v}\de\mathbf{x}\\
&\hspace{3cm} = \int_{\partial \Omega_N}{\bf g}\cdot{\bm v}\de\mathbf{x} 
\qquad\forall\,{\bm v} \in {\bm V}_0 \\
& \int_{\Omega}q\nabla \cdot{\bm u}(\mmu)\de\mathbf{x} =0 \qquad\forall\, q \in L^2(\Omega).  
\end{aligned}
\end{equation}
Problem \eqref{eq:weakNS-1} constitutes the particular case of the abstract problem \eqref{abs-pro} we use for the numerical illustrations.

\begin{figure}[h!]
\begin{center}
\includegraphics[width=4in]{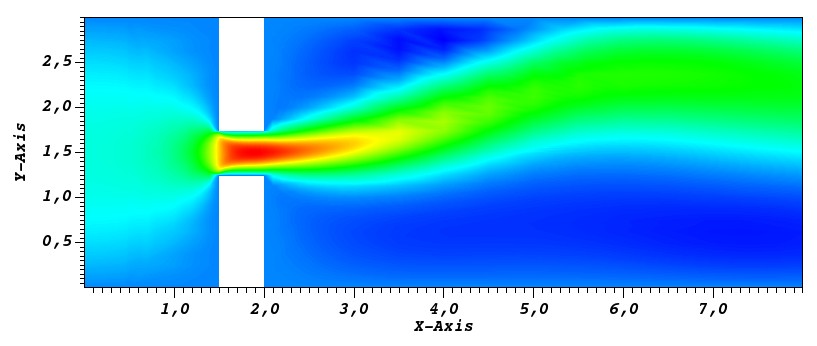} 
\includegraphics[width=.7in]{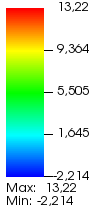} \\
\includegraphics[width=4in]{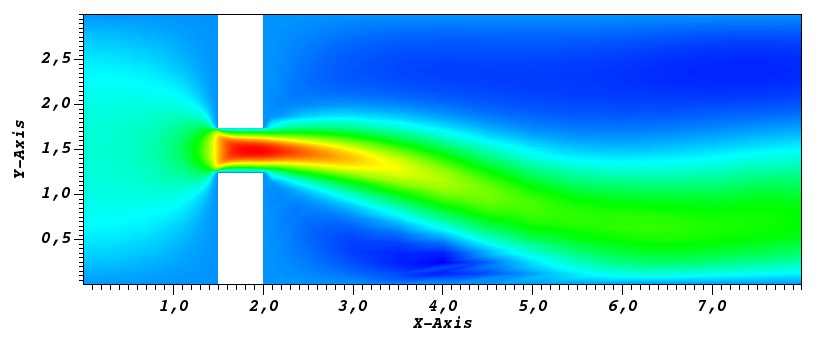}
\includegraphics[width=.7in]{figures/legend_max_para.png}
\caption{Steady state solutions for kinematic viscosity $\nu = 0.1$ and orifice width $w = 0.5$. 
Shown is the horizontal component of the velocity. 
The two stable solutions can be characterized by an attachment to the upper and lower wall, respectively.} 
\label{fig:channel_solution}
\end{center}
\end{figure}

We consider a benchmark test that has been widely studied in the literature:
channel flow through a narrowing of width $w$;
see, e.g., \cite{fearnm1,drikakis1,hawar1,mishraj1} and the references cited therein. 
The 2D geometry under consideration is depicted in Fig.~\ref{fig:channel_solution}. 
A parabolic horizontal velocity component with maximum $\frac{9}{4}$ and 
zero vertical component is inscribed on the inlet at the left side. 
At the top and bottom of the channel as well as the narrowing boundaries, zero velocity walls are assumed. 
The right end of the channel is an outlet, where zero Neumann boundaries \rev{(i.e., $\bf g = 0$)} are assumed.
We will let both the narrowing width and the viscosity vary in given ranges
that include a bifurcation. 

The Reynolds number $\text{Re}$ can be used to characterize the flow regime. 
For the chosen data, we have Re $= 9/ (4\nu)$\rev{, since the characteristic length has been set to one.}
As the Reynolds number $\text{Re}$ increases from zero, we first observe a steady {\em symmetric} jet with
two recirculation regions downstream of the narrowing that are symmetric about the centerline. 
As $\text{Re}$ increases, the recirculation length progressively increases. 
At a certain critical value $\text{Re}_{\text{crit}}$, one recirculation zone expands whereas the other shrinks, 
giving rise to a steady {\em asymmetric} jet. 
This asymmetric solution remains stable as $\text{Re}$ increase further, but the asymmetry becomes more pronounced. 
The configuration with a symmetric jet is still a solution, but is unstable \cite{sobeyd1}. 
Snapshots of the stable solutions for kinematic viscosity $\nu = 0.1$ and orifice width $w = 0.5$
are illustrated in  Fig.~\ref{fig:channel_solution}. 
This loss of symmetry in the steady solution as $\text{Re}$ changes is a supercritical pitchfork bifurcation \cite{Prodi}.

Because we are interested in studying a flow problem close to a steady bifurcation point, 
our snapshot sets include only steady-state solutions \cite{Hess2019CMAME}. 
To obtain the snapshots, we approximate the solution of problem \eqref{eq:weakNS-1} 
by a time-marching scheme that we stop when sufficiently close to the steady state, e.g., when the stopping condition

\begin{equation}
\frac{\|{\bm u}_N^n-{\bm u}_N^{n-1}\|_{L^2(\Omega)}}{\|{\bm u}_N^n\|_{L^2(\Omega)}}<\mathtt{tol}
\label{eq:def_increment}
\end{equation}

\noindent is satisfied for a prescribed tolerance $\mathtt{tol}>0$, where $n$ denotes the time-step index. 

\section{Numerical results}\label{sec:num_res}

We conduct a parametric study for the channel flow
where we let the the viscosity $\nu$ (physical parameter) vary in $[0.1, 0.2]$
and the narrowing width $w$ (geometric parameter) vary in $[ 0.5, 1.0 ]$.

\begin{figure}[h!]
\begin{center}
\includegraphics[width=.95\textwidth]{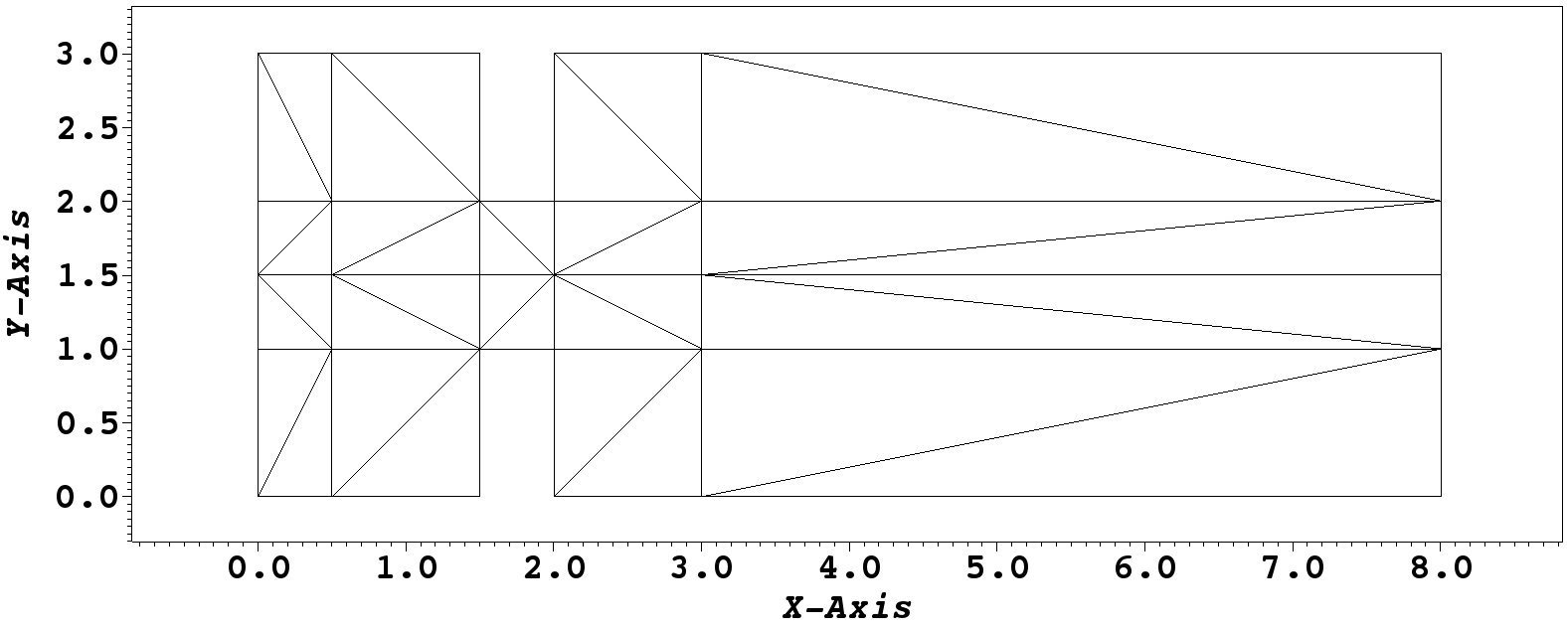}  
\caption{The $36$ triangular elements used for spatial approximation.}
\label{fig_SEM_domain_channel}
\end{center}
\end{figure}

We choose the Spectral Element Method (SEM) as FOM.
For the {\em spectral element discretization},
 the SEM software framework Nektar++, version 4.4.0, (see {\tt https://www.nektar.info/}) is used. 
The domain is discretized into $36$ triangular elements as shown in Fig. \ref{fig_SEM_domain_channel}. 
 Modal Legendre ansatz functions of order $12$ are used in every element and for every solution component. 
 This results in $4752$ degrees of freedom for each of the horizontal and vertical velocity components 
 and the pressure for the time-dependent simulations. For temporal discretization, an IMEX scheme of 
 order 2 is used with a time-step size of $\Delta t = 10^{-4}$; typically $10^5$ time steps are needed to reach a steady state.

\begin{figure}
\begin{center}
\includegraphics[width=1.0\textwidth]{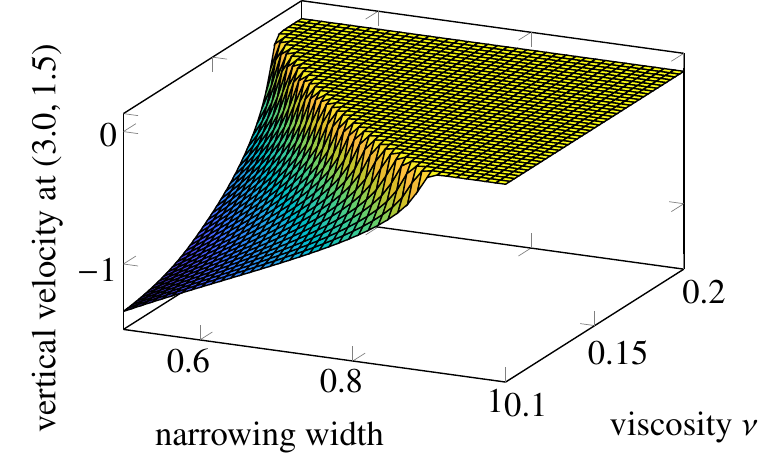}
\caption{Bifurcation diagram of the channel with variable width and viscosity:
 vertical component of the velocity computed by SEM close to steady state and evaluated at the point $(3.0, 1.5)$.
 Color encodes the vertical velocity component.}
\label{fig:bifurcation_diagram_channel_example_model_2p}
\end{center}
\end{figure}

Fig.~\ref{fig:bifurcation_diagram_channel_example_model_2p} shows
a bifurcation diagram: the vertical component of the velocity at the point $(3.0, 1.5) $ 
computed by SEM is plotted over the parameter domain.
The reference snapshots have been computed on a uniform $ 40 \times 41 $ grid.
Notice that Fig.~\ref{fig:bifurcation_diagram_channel_example_model_2p} reports only
the lower branch of asymmetric solutions. 

 In \cite{Hess2019CMAME}, we presented preliminary results obtained with our local ROM approach 
 for a two-parameter study related to the channel flow. 
 We showed that two criteria to select the local ROM basis that work well
 for one-parameter studies fail for the two-parameter case, in the sense
 that they provide a poor reconstruction of the bifurcation diagram. 
 Sec.~\ref{sec:criteria} addresses the need to find an accurate and inexpensive criterion to assign the
local ROM basis for a given parameter $\mmu$ in a multi-parameter context. 

The comparison for global ROM, local ROM, and POD-NN is reported in Sec.~\ref{sec:comparison}.

\subsection{Comparing criteria to select the local basis in the local ROM approach}\label{sec:criteria}

For our local ROM approach, we sample $72$ snapshots and divide them into 
8 clusters using k-means clustering. The number of clusters is chosen according to 
the minimal k-means energy. For more details, we refer to \cite{Hess2019CMAME}.

\begin{figure}
\begin{center}
\includegraphics[width=.9\textwidth]{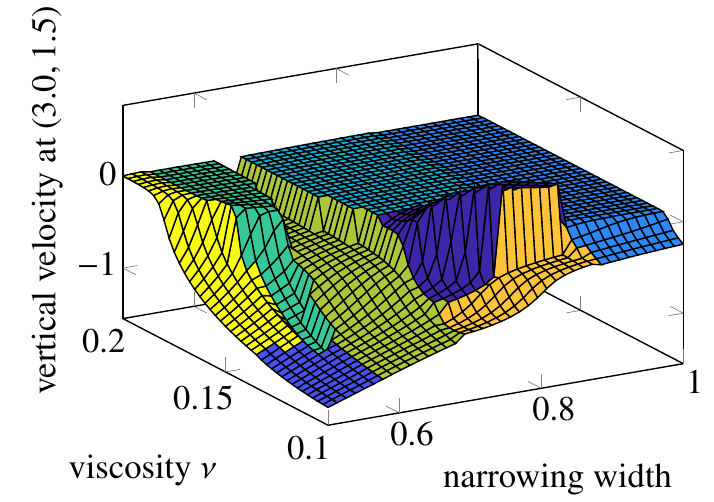}
\includegraphics[width=\textwidth]{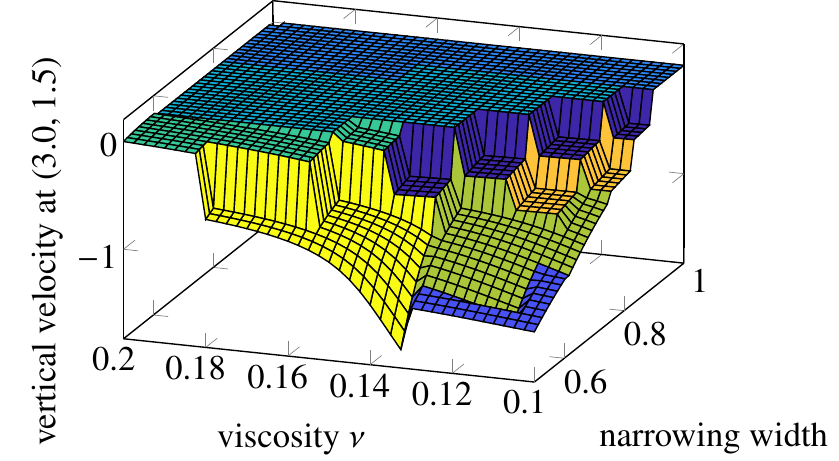}
\caption{Local ROM: bifurcation diagram reconstructed with
basis selection criterion distance to parameter centroid (top) and
distance to the closest snapshot location (bottom). 
Different colors are used for different clusters.}
\label{fig:compare_cluster_selection_cmame}
\end{center}
\end{figure}

\begin{figure}
\begin{center}
\includegraphics[width=\textwidth]{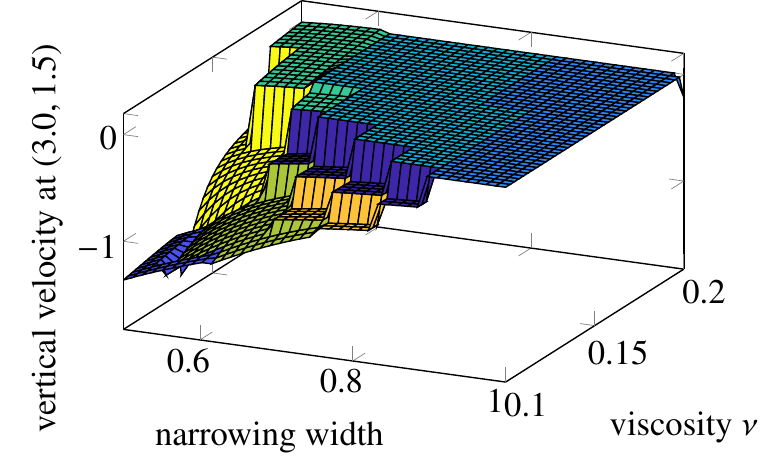}
\caption{Different view of the bifurcation diagram shown in 
Fig.~\ref{fig:compare_cluster_selection_cmame} (bottom).}
\label{fig:compare_cluster_no_overlap}
\end{center}
\end{figure}

First, we consider the two criteria presented in \cite{Hess2019CMAME}, namely
the  distance to parameter centroid and the distance to the closest snapshot location.
The first criterion entails finding the closest parameter centroid and 
using the corresponding local ROM basis.
The second criterion finds
the closest snapshot location to the given parameter vector $\mmu$ and the
local ROM basis that includes this snapshot is considered. 
The bifurcation diagrams reconstructed by the local ROM approach with
these two criteria are compared in Fig.~\ref{fig:compare_cluster_selection_cmame}.
The distance to parameter centroid criterion does not manage to recover the bifurcation diagram well 
and this seems largely due to the local ROM assignment scheme. See Fig.~\ref{fig:compare_cluster_selection_cmame} (top). 
In the bifurcation diagram corresponding to the distance to the closest snapshot location,
we observe several jumps when moving from one cluster to the next. 
See Fig.~\ref{fig:compare_cluster_selection_cmame} (bottom).
These jumps, which correspond to large approximation errors, can perhaps
be better appreciated from another view of the same bifurcation diagram shown
in Fig.~\ref{fig:compare_cluster_no_overlap}.
Next, we will try to reduce the jumps.

\begin{figure}
\begin{center}
\includegraphics[width=\textwidth]{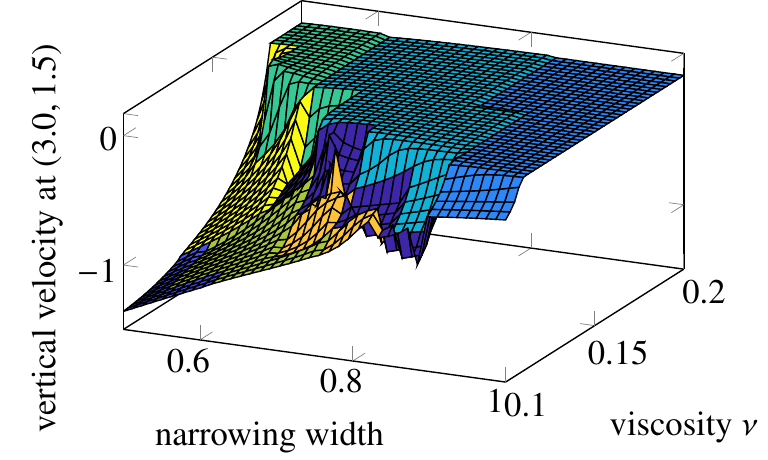}
\includegraphics[width=\textwidth]{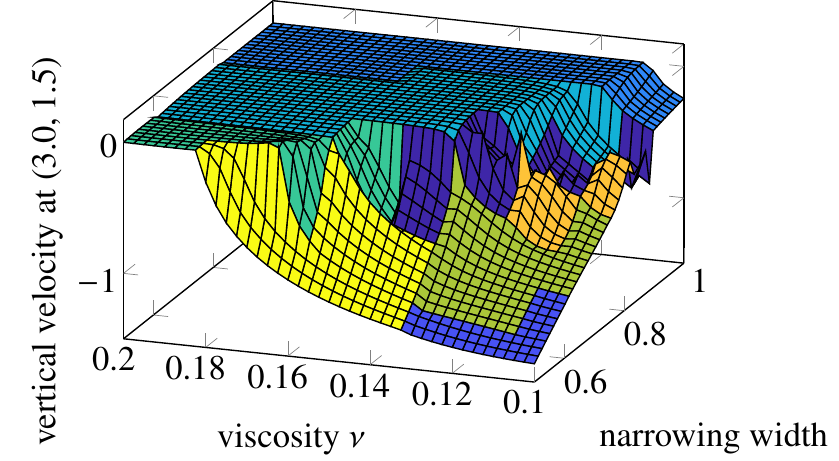}
\caption{Local ROM with basis selection criterion distance to the closest snapshot location and overlapping clusters:
two views of the reconstructed bifurcation diagram. Different colors are used for different clusters.}
\label{fig:compare_cluster_overlap}
\end{center}
\end{figure}

To alleviate the bad approximation in the transition regions between two clusters (i.e., local ROMs), 
we introduce overlapping clusters.
In particular, the collected snapshots of a cluster from the k-means algorithm undergo a first POD and 
are then enriched with the orthogonal complement of neighboring snapshots
according to the sampling grid. Then, a second POD with a lower POD tolerance  
defines the ROM ansatz space.
The corresponding bifurcation diagram is shown in Fig.~\ref{fig:compare_cluster_overlap}.
We observe that several jumps have been smoothed out. 
Compare Fig.~\ref{fig:compare_cluster_overlap} (top) with Fig.~\ref{fig:compare_cluster_no_overlap}.
The mean approximation error reduces by about an order of magnitude thanks to the overlap.

\begin{figure}
\begin{center}
\includegraphics[width=\textwidth]{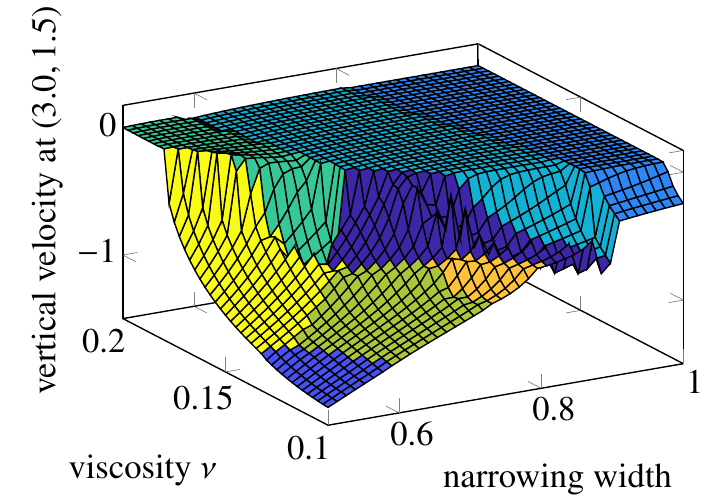}
\caption{Local ROM with basis selection selection that uses the ANN selection criterion: 
reconstruction of the bifurcation diagram.}
\label{fig:compare_cluster_selection_naive}
\end{center}
\end{figure}

Although the overlapping clusters lead to a better reconstruction of the bifurcation diagram, 
the result is still not satisfactory. Thus, we propose 
an alternative selection criterion that uses an artificial neural network (ANN). The ANN is trained using the k-means clustering as training information and enforcing a perfect match at the snapshot locations with the corresponding 
cluster. This is a classification problem, implemented in Keras Tensorflow.
\rev{The ANN is designed as multilayer perceptron with $4$ layers, the first layer having just two nodes (or "neurons") taking 
the two dimensional parameter values. The two inner layers are big ($2048$ and $1024$ nodes, respectively), 
while the last layer corresponds to the number of clusters, so $8$ in this case.
For the first three layers a ReLU activation function is used, while for the last layer activation function is a softmax.
The perfect match can be enforced by either running the training as long as the training data can be exactly matched or using a Keras early stopping with an outer loop checking for a match.
}
Fig.~\ref{fig:compare_cluster_selection_naive} shows the bifurcation diagram reconstructed with the
ANN selection criterion. We observe to a better reconstruction
of the bifurcation diagram: compare Fig.~\ref{fig:compare_cluster_selection_naive}
with Fig.~\ref{fig:compare_cluster_selection_cmame} and 
\ref{fig:compare_cluster_overlap}.


For a more quantitative comparison, we run 4 tests. The tests differ in the number of samples and
whether cluster overlapping is used. The specifications of each test are reported in Tables
\ref{table:errors_naive_criterion1} and \ref{table:errors_naive_criterion2}, which 
list the relative errors evaluated over the $40 \times 41$ reference grid 
in the $L^2$ and $L^\infty$ norms, respectively. Three local basis selection 
criteria are considered: distance to parameter centroid, distance to the closest snapshot location,
and ANN.
Tables \ref{table:errors_naive_criterion1} and \ref{table:errors_naive_criterion2} confirm that the overlapping
cluster represent an improvement over non-overlapping clusters. This is true for all three criteria, but in particular
for the distance to parameter centroid criterion. Thus, for tests 3 and 4 we only used overlapping clusters.
We notice that the ANN criterion outperforms the other two criteria in all the tests, both in the 
$L^2$ and $L^\infty$ norms. However, 
the margin of improvement becomes smaller as the number of samples increases. 
Tables \ref{table:errors_naive_criterion1} and \ref{table:errors_naive_criterion2} also report
the mean relative error for an optimal cluster selection, which is explained next.

\begin{table}[h!]
\begin{center}
    \begin{tabular}{ | l | l | l | l | l |}
    \hline
     & Test 1 & Test 2 & Test 3 & Test 4 \\ \hline
    samples & 72 & 72 & 110 & 240 \\ \hline
    uniform grid & $8 \times 9$ & $8 \times 9$ & $10 \times 11$ & $15 \times 16$ \\ \hline
    overlapping clusters & yes & no & yes & yes \\ \hline
    parameter centroid mean & 0.0294 & 0.0632 & 0.1085 & 0.0559 \\ \hline
    distance snapshot mean & 0.0241 & 0.0295 & 0.1011 & 0.0227 \\ \hline
    ANN mean  & 0.0238 & 0.0273 & 0.1002 & 0.0223 \\ \hline
    optimum  & 0.0046 & 0.0106 & 0.0048 & 0.0092 \\ \hline
    \end{tabular}
    \caption{Local ROM with three different basis selection: mean 
    relative $L^2$ errors for the velocity over all reference parameter points ($1640$ on a uniform $40 \times 41$ grid).}
   \label{table:errors_naive_criterion1}
   \end{center}
\end{table}
\begin{table}[h!]
\begin{center}
    \begin{tabular}{ | l | l | l | l | l |}
    \hline
     & Test 1 & Test 2 & Test 3 & Test 4 \\ \hline
    samples & 72 & 72 & 110 & 240 \\ \hline
    uniform grid & $8 \times 9$ & $8 \times 9$ & $10 \times 11$ & $15 \times 16$ \\ \hline
    overlapping clusters & yes & no & yes & yes \\ \hline
    parameter centroid mean  & 0.0275 & 0.0600 & 0.0970 & 0.0538 \\ \hline
    distance snapshot mean & 0.0230 & 0.0284 & 0.0908 & 0.0217 \\ \hline
    ANN mean & 0.0227 & 0.0263 & 0.0899 & 0.0214 \\ \hline
    optimum  & 0.0044 & 0.0101 & 0.0045 & 0.0088 \\ \hline
    \end{tabular}
    \caption{Local ROM with three different basis selection: mean relative $L^\infty$ errors for the velocity
    over all reference parameter points ($1640$ on a uniform $40 \times 41$ grid).}
   \label{table:errors_naive_criterion2}
   \end{center}
\end{table}

\begin{figure}
\begin{center}
\includegraphics[width=.75\textwidth]{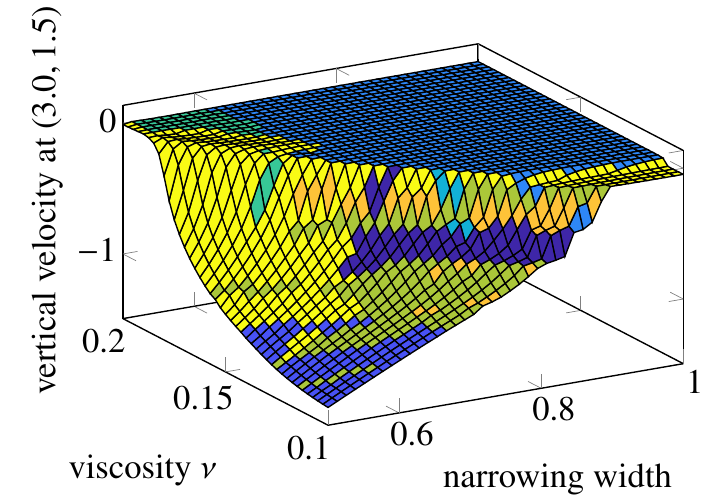}
\caption{Optimal cluster selection. Once again, different colors are used for different clusters.}
\label{fig:cluster_selection_optimal}
\end{center}
\end{figure}

\rev{By "optimal cluster selection" we address the question of how parameter points are optimally associated with an already given clustering.}
Fig.~\ref{fig:cluster_selection_optimal} shows the optimal cluster selection. 
This optimal selection can only be obtained from a fine grid of reference solutions. Thus it is usually not available. 
The reason why we report it is because some interesting conclusions can be drawn from it. 
First, the best possible clustering does not have contiguous clusters.
This is in contrast to the clusters created by the k-means algorithm, which are contiguous in the parameter space in all of our tests.
Second, the best possible approximation at a snapshot location is not necessarily given by the local cluster 
to which that snapshot is assigned.
Third, there still is a reduction factor of $5 - 20$ in the relative $L^2$ for the velocity error that 
could be gained, as one can see when comparing the respective errors 
for the optimal cluster and the ANN selection criteria 
in  Tables~\ref{table:errors_naive_criterion1} and \ref{table:errors_naive_criterion2}.


To get close to the optimum, we adopt the following strategy.  
We compute relative errors of all local ROMs at all snapshot locations. 
This operation is performed offline and is not expensive since the exact 
solution at the snapshot locations is available.
The relative errors can be used as training data for an ANN, which means that 
the ANN training is treated as a regression and not a classification. Thus,
the ANN will approximate relative errors of each local ROM over the parameter domain.
This approximation is used as cluster selection criterion.
In this procedure, it is important to normalize the data.  
Here, we use the inverse relative error and  normalize the error vectors at each snapshot location.
We note that the training time for the ANN increases, but is still well below the ROM offline time, i.e.~30 minutes vs. several hours. 
\rev{The ROM offline time is dominated by computing the affine expansion of the trilinear form. The computational cost grows 
with the cube of the reduced order model dimension. This makes the localized ROM much faster than the global ROM as, 
for example, the global ROM might have dimension of $40$ while each of the
$8$ local ROMs has a dimension of about $10$. We found impossible to quantify
the training cost of an ANN. The performance of the stochastic gradient employed in the ANN training
varied significantly over multiple runs and required outer loops to check the accuracy. Thus, the given 
measure of ~30 minutes vs. several hours is only our experience with this particular model and probably cannot be generalized.}

We consider test 3 ($10 \times 11$ sampling grid, overlapping clusters) to assess two variants of the regression ANN, 
which differ in how the training is done. 
One variant takes all clusters into account simultaneously and is called ``regression ANN''. 
It generates a mapping from $\mathbb{R}^2 \mapsto \mathbb{R}^8$, i.e.~the two dimensional 
parameter domain to the expected errors of the clusters.
The second variant treats each clusters independently and is called ``regression ANN, independent local ROMs''. 
It generates eight mappings from $\mathbb{R}^2 \mapsto \mathbb{R}$, i.e.~one for each cluster.
Table~\ref{table:errors_smart_crit} reports the mean relative $L^2$ and $L^\infty$ errors for the velocity. 
Interestingly, taking all clusters into account simultaneously (``regression ANN'') is about $10\%$ 
more accurate than considering each cluster separately (``regression ANN, independent local ROMs''). 
Moreover,  we observe that the mapping snapshot location to errors of local ROMs 
holds useful information and the ANN consistently gets closest to the optimum.
Table~\ref{table:errors_smart_crit} reports also the errors obtained with the
Kriging DACE software\footnote{Kriging DACE - Design and Analysis of Computer Experiments, \texttt{http://www.omicron.dk/dace.html}.}.
and with simply taking the cluster, which best approximates the closest snapshot. 
The closest snapshot is determined in parameter domain with the Euclidean norm. 
This is an easily implementable tool, which still performs better than the distance to the parameter centroid; 
see  \cite{Hess2019CMAME}. 
From Table~\ref{table:errors_smart_crit}, we see that this simple criterion performs 
only $15\% - 20\%$ worse and is very cheap to evaluate.

  \begin{center}
\begin{table}[h!]
    \begin{tabular}{ | l | l | l | }
    \hline
     & mean $L^2$ error & mean $L^\infty$ error  \\ \hline
    optimum & 0.0048 & 0.0045  \\ \hline
    regression ANN & 0.0068 & 0.0064  \\ \hline
    regression ANN, independent local ROMs & 0.0076 & 0.0071  \\ \hline
    Kriging DACE & 0.0077 & 0.0073  \\ \hline
    distance to next best-approx. snapshot  & 0.0081 & 0.0077  \\ \hline
    \end{tabular}
    \caption{Local ROM with different basis selection: mean relative $L^2$ and $L^\infty$
    errors for the velocity over all reference parameter points ($1640$ on a uniform $40 \times 41$ grid).}
    \label{table:errors_smart_crit}
\end{table}
\end{center}

\subsection{Comparing the global ROM, local ROM, and POD-NN approaches}\label{sec:comparison}

In the previous section, we learned that the regression ANN criterion outperforms all other local basis 
selection criteria. In this section, we compare the local ROM approach with the regression ANN criterion
to our global ROM approach and the POD-NN over the reconstruction of the bifurcation diagram. 

Table \ref{table:errors_compare_to_pod_nn} reports the mean relative $L^2$ and $L^\infty$ errors
for the velocity for the three approaches under consideration. We consider four different
numbers of samples: 42, 72, 110, and 240.
It can be observed that the global ROM shows a slow convergence, not even $7\%$ accuracy is reached with 
the finest sampling grid. 
The local ROM with regression ANN cluster selection shows no distinctive convergence behavior, 
which might indicate that the accuracy saturates at lower snapshot grid sizes. 
The POD-NN shows the fastest convergence, reaching about $0.3\%$ error with the finest sampling grid. 
We note that the POD-NN training did not take overfitting into account.
\rev{Overfitting occurs when the training data is more accurately approximated than the actual data of interest. 
It can be checked with having a validation set, whose accuracy is measured independently and not included in the training.
Training can be stopped when the training data is more accurately approximated than the validation set.}

\begin{table}[h!]
  \begin{center}
    \begin{tabular}{ | l | l | l | }
    \hline
    42 snapshots  & mean $L^2$ error & mean $L^\infty$ error  \\ \hline
    global ROM  & 3.7022 & 3.1120  \\ \hline
    local ROM + regression ANN & 0.0510 & 0.0486  \\ \hline 
    POD-NN  & 0.0108 & 0.0104  \\ \hline   \hline
    72 snapshots & mean $L^2$ error & mean $L^\infty$ error  \\ \hline
    global ROM& 0.6970 & 0.5831  \\ \hline
    local ROM + regression ANN  & 0.0103 & 0.0098  \\ \hline 
    POD-NN & 0.0080 & 0.0075  \\ \hline \hline
    110 snapshots   & mean $L^2$ error & mean $L^\infty$ error  \\ \hline
    global ROM & 0.1044 & 0.0948  \\ \hline
    local ROM + regression ANN & 0.0068 & 0.0064  \\ \hline
    POD-NN & 0.0059 & 0.0053  \\ \hline \hline
    240 snapshots   & mean $L^2$ error & mean $L^\infty$ error  \\ \hline
    global ROM & 0.0762 & 0.0734  \\ \hline
    local ROM + regression  ANN & 0.0101 & 0.0096  \\ \hline
    POD-NN & 0.0032 & 0.0027  \\ \hline
    \end{tabular}
    \end{center}
    \caption{Comparison of global ROM, local ROM and POD-NN for four snapshot grids.}
    \label{table:errors_compare_to_pod_nn}
\end{table}

Several remarks are in order.

\begin{remark}
The data reported in Table \ref{table:errors_compare_to_pod_nn} concerning the local ROM with regression ANN cluster selection
and the POD-NN are sensitive to the neural network training. In addition, the data for local ROM approach
are sensitive to the POD tolerances are changed. Nonetheless, Table \ref{table:errors_compare_to_pod_nn}
provides a general indication of the performance of each method in relation to the other two.  
\end{remark}

\begin{remark}
A rigorous comparison in term of computational times is not possible because the different
methods are implemented in different platforms. However, we can make some general comments. 
The POD-NN method has a significant advantage in terms of computational time:
one does not need to assemble the trilinear reduced form associated to the convective term, which 
is our simulations takes about 1-3 hours. The time required for the
POD-NN evaluation in the online phase is virtually zero, while
projection methods need to do a few iterations of the reduced fixed point scheme. That takes about 10-30 s.
On the other hand, the POD-NN required the training of the ANN. However, in our simulations
that takes only about 20 minutes. 
\end{remark}

\begin{remark}
We also investigated a local POD-NN, i.e. we combined a k-means based localization approach with the POD-NN. 
However, this led to significantly larger errors than the global POD-NN method.
In general, the error of a local POD-NN approach was in the range of the error of the ``local ROM + regression ANN''. Thus, 
we did not pursue this approach any further. 
\end{remark}

\section{Concluding remarks}\label{sec:concl}
We focused on reduced-order models (ROMs) for PDE 
problems that exhibit a bifurcation when more than one parameter is varied.
For a particular fluid problem that features a supercritical pitchfork bifurcation 
under variation of Reynolds number and geometry,
we investigated projection-based local ROMs and compared them
to an established global projection-based ROM as well as an emerging artificial neural network (ANN) based method
called POD-NN. 
We showed that k-means based clustering, transition regions, cluster-selection criteria based on best-approximating clusters and 
ANNs gain more than an order of magnitude in accuracy over the global projection-based ROM. 
Upon examining the accuracy of POD-NN, it became obvious that the POD-NN provides consistently 
more accurate approximations than the local projection-based ROM. 
Nevertheless, the local projection-based ROM might be more amenable to the use of reduced-basis error estimators than the POD-NN.
This could be the object of future work.

\bibliographystyle{plain}
\bibliography{rbsissa,latexbi}

\end{document}